\documentclass[a4paper, oneside, 11pt]{amsart}
\usepackage[utf8]{inputenc}
\usepackage[english]{babel}

\usepackage[big]{layaureo}

\usepackage{amssymb, latexsym}

\usepackage[colorlinks=true, linkcolor=red, urlcolor=black, citecolor=blue]{hyperref}

\DeclareMathOperator{\id}{\ensuremath{id}}

\newcommand{\HH}{\ensuremath{\mathbb{H}}} 
\newcommand{\C}{\ensuremath{\mathbb{C}}} 
\newcommand{\R}{\ensuremath{\mathbb{R}}} 
\newcommand{\Z}{\ensuremath{\mathbb{Z}}} 
\newcommand{\calE}{\ensuremath{\mathcal{E}}} 
\newcommand{\calH}{\ensuremath{\mathcal{H}}} 
\newcommand{\calQ}{\ensuremath{\mathcal{Q}}} 
\newcommand{\calU}{\ensuremath{\mathcal{U}}} 
\newcommand{\Sph}{\mathrm{S}}
\newcommand{\Ker}{\mathrm{Ker}}

\newcommand{\tr}{\mathrm{tr}}

\newcommand{\Bis}{\mathrm{B}}	
\newcommand{\SU}{\mathrm{SU}}
\newcommand{\U}{\mathrm{U}}
\newcommand{\Sp}{\mathrm{Sp}}

\newcommand{\scal}{\mathrm{s}}
\newcommand{\can}{\mathrm{can}}
\newcommand{\frakd}{\mathfrak{d}}		
\newcommand{\frakf}{\mathfrak{f}}

\newcommand{\fraksu}{\mathfrak{su}}
\newcommand{\RiC}{\mathrm{\rho}}	
\newcommand{\eps}{\epsilon}

\newcommand{\img}{\mathtt{i}}
\newcommand*{\dd}{\mathop{}\!\mathrm{d}} 			
\newcommand*{\del}{\mathop{}\!\mathrm{\partial}} 	
\newcommand{\hc}{hyper-complex}		
\newcommand{\hh}{hyper-Hermitian}	
\newcommand{\hk}{hyper-Kähler}		
\newcommand{\hn}{hyper-normal}		
\renewcommand{\Re}{\mathrm{Re}}

\renewcommand{\phi}{\varphi}

\providecommand{\abs}[1]{\lvert#1\rvert}

\theoremstyle{plain}
\newtheorem{thm}{Theorem}[section]
\newtheorem{prop}[thm]{Proposition}
\newtheorem{lem}[thm]{Lemma}
\newtheorem{cor}[thm]{Corollary}

\theoremstyle{definition}

\newtheorem{rmk}[thm]{Remark}
\newtheorem{ex}[thm]{Example}

\title{HKT-Einstein metrics on circle bundles}
\author{Lorenzo Marcocci}
\date{}
\thanks{This work was partially supported by GNSAGA of INdAM}
\address{Dipartimento di Ingegneria e Scienze dell'Informazione e Matematica, Università dell'Aquila, Via Vetoio, 67100 L'Aquila, Italy}
\email{lorenzo.marcocci@graduate.univaq.it}

\begin{document}

\maketitle

\begin{abstract}
We consider a construction of \hh{} structures on total spaces of circle bundles over 3-quasi-Sasakian manifolds.
We study the HKT-Einstein equation in this setting, proving that a large family of non-homogeneous \hc{} manifolds admit compatible HKT-Einstein metrics with positive Einstein factor.
We also show that the \hh{} structures on such circle bundles are quaternionic balanced or quaternionic Gauduchon if and only if they are HKT.
Finally, the HKT-Einstein condition is investigated under the twist construction.
\end{abstract}

\section{Introduction}\label{S:Intro}
A \emph{\hc{} structure} on a smooth manifold $M$ is a triple of complex structures $(J_1, J_2, J_3)$ obeying the relations of unit quaternions: 
$J_1 J_2 = - J_2 J_1 = J_3$.

Clearly, \emph{\hk{} manifolds} belong to this family.
However, the restrictions imposed by the \hk{} condition led to the study of more general notions of special \hh{} metrics.
Among them, a prominent role is played by \emph{HKT manifolds}, which firstly appeared in relation to theoretical physics (\cite{HP}) and are still very important in that field of study.
They are, in many ways, the \hc{} counterpart of Kähler manifolds.
More recently, in \cite{FusiGentili} the authors proposed an Einstein-type condition in this setting, introducing the notion of \emph{HKT-Einstein metrics}.
Under the above analogy, they may correspond to Kähler-Einstein metrics.

At present, the class of known examples of HKT-Einstein metrics with positive Einstein factor is still rather small.
There exist invariant HKT-Einstein metrics on compact  homogeneous \hc{} manifolds: we refer the reader to~\cite{FusiGentili} for the case of Lie groups, and~\cite{BM} for the general case and for uniqueness.
Non-homogeneous examples appeared recently on the quaternionic Hopf surface, or as deformations of locally conformally \hk{} metrics (\cite{BrienzaStreets, GentiliVezzoni}).\\

In this work, we consider a specific class of \hh{} manifolds, obtained as \emph{total spaces of principal circle bundles over 3-quasi-Sasakian manifolds}.
The construction of Hermitian manifolds on circle bundles over almost-contact manifolds goes back to Ogawa~\cite{Ogawa}.
In the \hc{} case, it was first considered  in~\cite{Cabrera, FFUV1, FFUV2} and obtained independently in \cite{BGM-Stiefel, BGM-hc}.
Each of these papers focuses on different aspects.
In particular, the property of being (strong) HKT was characterized  in~\cite{FFUV1}.

Our first purpose is to find examples of HKT-Einstein metrics, by studying the HKT-Einstein equation on these manifolds.
This project is partly motivated by the fruitful research of Einstein metrics on circle bundles (see e.g.~\cite{WangZiller}).
The central results are contained in Section~\ref{S:HKTE}. In particular, in Subsection~\ref{SS:inhom} we apply our analysis to reveal the existence of an infinite family of compact simply connected HKT-Einstein manifolds, which are not homogeneous.
The underlying \hc{} manifolds were introduced in~\cite{BGM-Stiefel}, 
and the bases of the fibrations have been obtained via the method of 3-Sasakian reduction (\cite{BGM-Main}).

The same approach also produces HKT manifolds that are not HKT-Einstein, but have constant Chern scalar curvature.\\

Besides HKT structures, there are other special classes of \hh{} metrics, namely \emph{quaternionic Gauduchon} and \emph{quaternionic balanced} metrics; they are important in their own right for the geometry of \hc{} manifolds (see~\cite{FusiGentili, GLV}).
Since the \hh{} structures on circle bundles considered here are not always HKT, it is natural to ask if, in those cases, they satisfy these weaker conditions.
In Section~\ref{S:qBal},  we give a negative answer to this question, proving that the metric can be quaternionic Gauduchon or quaternionic balanced if and only if it is HKT.\\

In Section~\ref{S:Twist} we go back to the HKT-Einstein condition, with the aim to study it in relation to the twist construction introduced by A. Swann in~\cite{Swann}.
This technique has proven to be valuable in several cases, in complex and \hc{} geometry.
Here we show that, with an assumption on the twist data, the twist of a HKT-Einstein manifold has again this property.
To produce new examples, we then show how the twist can be applied to the product between one of the manifolds studied in Section~\ref{S:HKTE}, and a homogeneous HKT-Einstein manifold.\\

\textbf{Acknowledgements.}
The author thanks Lucio Bedulli and Giovanni Gentili for reading an earlier version of the manuscript and providing useful comments and suggestions.

\section{Preliminaries}\label{S:Prelim}
\subsection{Hyper-Hermitian geometry}\label{SS:HHerm}
Let $(M^{4n}, J_1, J_2, J_3)$ be a \hc{} manifold.
A Riemannian metric $h$ on $M$ is called \hh{} provided it is Hermitian with respect to each complex structure $J_i$.
Such a  metric  induces a corresponding form
\[
	\Omega := \frac{1}{2} \left(\omega_{J_2} + \img\omega_{J_3} \right),
\]
which lies in $\Lambda^{2,0}_{J_1}(M)$ .
Here $\omega_{J_i} = h(J_i \cdot, \cdot)$.
One can recover $h$ from $\Omega$ via the relation $h(X,Y) = 2 \Re (\Omega(X, J_2Y))$.

An important family  is given by \emph{\hk{} metrics with torsion}, in short \emph{HKT metrics}, which are characterized by the equality of  the Bismut connections of $(J_i, h)$, $i=1,2,3$.
In other words, $J_1\dd\omega_{J_1} = J_2\dd\omega_{J_2} = J_3\dd\omega_{J_3}$.
As proved in~\cite{GrantcharovPoon}, this condition is equivalent to
\[
	\del\Omega = 0,
\]
where the splitting $\dd = \del + \bar{\del}$ is taken with respect to $J_1$.
\emph{Strong HKT metrics} are defined requiring the torsion of the Bismut connection to be closed, that is, each Hermitian structure has to be SKT.

Also other special metrics appeared in the context of \hc{} geometry.
A \hh{} metric $h$ is called
\begin{itemize}
\item
\emph{quaternionic balanced} if $\del\Omega^{n-1} = 0$;
\item
\emph{quaternionic Gauduchon} if $\del{\del_{J_2}\Omega^{n-1}} = 0$, where $\del_{J_2} = J_2^{-1}\bar{\del}J_2$.
\end{itemize}
These notions have been introduced respectively in~\cite{LejmiWeber} and~\cite{GLV}, as analogues of \emph{balanced} and \emph{Gauduchon} metrics in this  setting.
Note that quaternionic Gauduchon is weaker than quaternionic balanced, and both are weaker than  HKT.

Following~\cite{FusiGentili}, a HKT metric $h$ is called \emph{HKT-Einstein} when there exists $\lambda\in\R$ such that
\begin{equation}\label{E:HKT-E}
	\frac{\RiC_{J_1} - J_2 \RiC_{J_1} }{2} = \lambda \omega_{J_1}.
\end{equation}
Here and below $\RiC_{J_i}$ is the first Chern--Ricci form associated to $(J_i, h)$, defined by
\[
	\RiC_{J_i}(X,Y) 
	= \sum_{k=1}^{2n} 
		h\!\left( R(X,Y) J_ie_k, e_k \right)
\]
where $\{ e_k, J_ie_k \}_{k=1}^{2n}$ is a local orthonormal frame, and $R$ denotes the curvature of the Chern connection.
It is known that on compact manifolds necessarily $\lambda \geq 0$, 
and moreover $\lambda = 0$ if and only if $h$ is balanced with respect to one (hence all) of the complex structures $J_i$.

Here are other properties of \hh{} manifolds that we shall use; proofs of them can be found in~\cite{FusiGentili}.
\begin{enumerate}
\item
There exist $ \beta_{\Omega} \in \Lambda_{J_1}^{1,0}(M)$ such that $\del\Omega^{n-1} = \beta_{\Omega}\wedge\Omega^{n-1}$.
The quaternionic balanced condition is equivalent to $\beta_{\Omega}=0$, and it holds that
\[
	\beta_{\Omega} = \Lambda(\del\Omega),
\]
where $\Lambda$ is the adjoint of the Lefschetz operator $L= \Omega\wedge(-)$.
\item
The Lee forms of the individual Hermitian structures all coincide with the same 1-form $\tau$, and it holds that $\rho_{J_i}=\dd J_i( \tau - \beta_{\Omega} - \bar{\beta}_{\Omega})$. In particular, if the manifold is HKT then 
\[
	\rho_{J_i} = \dd J_i\tau.
\]	
\item
The Chern scalar curvatures of the individual Hermitian structures all coincide with the same function, $\scal$. The same holds for the Bismut scalar curvature, that we denote by $\scal^{\Bis}$.
\item
For a HKT-Einstein manifold with Einstein factor $\lambda$, it holds that
\[
	\scal = 2n\lambda,
\]
where $4n$ is the real dimension. Therefore, HKT-Einstein manifolds fall in the larger class of HKT manifolds with \emph{constant Chern scalar curvature}.
\end{enumerate}

\subsection{Almost 3-contact manifolds}\label{SS:a3cm}

An \emph{almost contact} structure on a manifold $N^{2m+1}$ is a triple consisting of  a $(1,1)$-type tensor field $\phi$, a vector field $\xi$, called \emph{Reeb vector field}, and a 1-form $\eta$ satisfying
\[
	\phi^2 = - \id + \eta\otimes\xi, \quad \eta(\xi)=1.
\]
These conditions imply $\phi(\xi)=0$ and $\eta\circ\phi=0$.
The tangent bundle $TN$ decomposes as $TN = \mathcal{D} \oplus \langle \xi \rangle$, where $\mathcal{D} = \Ker\,\eta = \operatorname{Im} \phi$.
The product manifold $N\times\R$ can be endowed with an almost complex structure $\bar{J}$ given by
\[
	\bar{J} \left( X, f \frac{\dd}{\dd t} \right)
	= \left( \phi X + f \xi,  - \eta(X) \frac{\dd}{\dd t} \right),
\]
where $f$ is a smooth function on $N\times\R$.
The structure is called \emph{normal} if $\bar{J}$ integrable, and by~\cite{SH} this is equivalent to
\[
	[\phi, \phi] + \dd\eta\otimes\xi = 0,
\]
where $[\phi, \phi]$ is the Nijenhuis tensor of $\phi$, given by $[\phi, \phi](X,Y) = \phi^2 [X,Y]+[\phi X, \phi Y]- \phi[\phi X,Y]- \phi[X, \phi Y]$.
It is well known that for a normal almost contact structure we have $\iota_{\xi}\dd\eta = 0$, see e.g.~\cite{Blair}.

An almost contact structure is called \emph{metric} if there is additionally a Riemannian metric $g$ such that
\[
	g(\phi X, \phi Y) = g(X,Y) - \eta(X)\eta(Y)
	\quad
	\forall \, X,Y \in TN.
\]
Note that the sum $TN = \mathcal{D}\oplus\langle \xi \rangle$ is orthogonal, $\xi$ has unit length, and $\eta$ is dual to $\xi$ with respect to $g$.
The \emph{fundamental form} $\Phi$ is defined as $\Phi(X,Y) = g(X, \phi Y)$.

Important classes of normal almost contact metric manifolds are given by \emph{quasi-Sasakian} manifolds, defined by the condition $\dd\Phi=0$,
and \emph{$\alpha$-Sasakian} manifolds, which verify $\dd\eta = 2\alpha\Phi$ for some $\alpha\in\R^*$. \\

An \emph{almost 3-contact structure} on a smooth manifold $N^{4m+3}$ is given by three almost contact structures $(\phi_i, \xi_i, \eta_i)$, $i = 1,2,3$, such that
\[
\begin{gathered}
	\phi_i \circ \phi_j - \eta_j \otimes \xi_i = - \phi_j \circ \phi_i + \eta_i \otimes \xi_j 
		= \eps_{ijk}\phi_k, \\
	\phi_i(\xi_j) = - \phi_j(\xi_i) = \eps_{ijk}\xi_k,\\
	\eta_i \circ \phi_j = - \eta_j \circ \phi_i = \eps_{ijk}\eta_k,
\end{gathered}
\]
for every distinct $i,j,k$, where $\eps_{ijk}$  is the Levi-Civita symbol. 
If each almost contact structure is normal, the structure is said to be \emph{\hn}.
The tangent bundle decomposes as $TN = \calQ \oplus \calU$, where
\[
	\calQ \, :=\, \bigcap_{i=1}^{3}\operatorname{Ker}\eta_i,\qquad
	\calU \, :=\, \langle\xi_1,\xi_2,\xi_3\rangle.
\]
Any almost $3$-contact manifold admits a Riemannian metric $g$ which is compatible with
each of the three structures, and in this case one refers to a \emph{metric} almost $3$-contact manifold.
The subbundles $\calQ$ and $\calU$ are $g$-orthogonal and the Reeb vector fields are orthonormal.
From now on, we shall always assume that $m\geq1$, so that $\calQ$ has rank at least one.\\

We will be mainly interested in \emph{3-quasi-Sasakian manifolds}, namely almost 3-contact metric manifolds such that each structure $(\phi_i, \xi_i, \eta_i, g)$ is quasi-Sasakian.
A number of  properties of these manifolds have been established  in~\cite{CMDND1, CMDND2}, and we refer the reader to these papers for detailed proofs.
First of all, there is a constant $\delta\in\R$ such that 
\begin{equation}\label{E:delta}
	[\xi_i, \xi_j] = 2\delta\eps_{ijk}\xi_k \quad \text{for every $i,j,k$.}
\end{equation}
There exists an odd integer $2p + 1$, called \emph{rank} of the manifold, such that $\eta_i \wedge (\dd\eta_i)^p \neq 0$ and $(\dd\eta_i)^{p+1} = 0$, for every $i = 1,2,3$.
Moreover the rank is of the form $4l+3$ or $4l+1$, $0\leq l \leq m$, according to $\delta\neq0$ or $\delta=0$.

The distribution $\calQ$ decomposes as $\calQ = \calE^{4l}\oplus\calE^{4q}$, where $\calE^{4q} := \{ X \in \calQ \colon \iota_X\dd\eta_i = 0 \text{ for every $i$} \}$ and $\calE^{4l}$ is its orthogonal complement in $\calQ$ (here $l+q=m$).
It holds that $\phi_i(\calE^{4l}) = \calE^{4l}$ for every $i$, and on $\calE^{4l}$ 
we have $\dd\eta_i = 2\delta\Phi_i$ for each $i$.

A 3-quasi-Sasakian manifold is called
\begin{itemize}
\item
\emph{3-$\delta$-cosymplectic} if $\dd\eta_i = -2\delta \eps_{ijk}\eta_j \wedge \eta_k$ for every $i,j,k$,
and \emph{3-cosymplectic} if $\delta = 0$;
\item
\emph{3-$\alpha$-Sasakian} if $\dd\eta_i = 2\alpha\Phi_i$ for every $i$, where $\alpha\in\R^*$,
and \emph{3-Sasakian} if $\alpha = 1$. 
\end{itemize}
In fact, 3-quasi-Sasakian manifold is 3-cosymplectic if and only if the constant $\delta$ given by~\eqref{E:delta} vanishes.
We should also emphasize that every 3-quasi-Sasakian manifold with $\delta\neq0$ is locally the Riemannian product of a 3-$\alpha$-Sasakian manifold, tangent to the distribution $\calE^{4l}\oplus\mathcal{U}$, and a \hk{} manifold, tangent to $\calE^{4q}$  (see in particular~\cite{CMDND2}).\\

Finally, for a differential $p$-form $\gamma$ on $N$, we set
\[
	(\phi_i \cdot \gamma)(X_1, \dots, X_p) := \gamma(\phi_i X_1, \dots, \phi_i X_p).
\]
It holds that, for every $X,Y\in\Gamma(\calQ)$ and every even permutation $(ijk)$ of $(123)$:
\begin{equation}\label{E:formulas}
	\iota_{\xi_j}\dd\eta_i = -2\delta\eta_k,\qquad
	\dd\eta_i(\phi_j X,Y) = \dd\eta_k(X,Y),\qquad
	\phi_i \cdot \dd\eta_j = - \dd\eta_j + 2\delta \eta_i \wedge \eta_k.
\end{equation}

\subsection{Hyper-Hermitian structures on circle bundles}\label{SS:cb}
The following construction is originally due to Ogawa \cite{Ogawa}.
Consider  a smooth manifold $N^{2m+1}$ with an almost contact structure $(\phi,\xi,\eta)$.
Let $\pi\colon M^{2n} \to N^{2m+1}$ be a principal $\Sph^1$-bundle over $N$,  where $\Sph^1 = \R/\Z$, and let $\theta$ be a connection 1-form on it. 
Since the structural group of the bundle is $\Sph^1$, the curvature form will be $\dd\theta = \pi^{*} F$, for a 2-form $F$ on $N$.

It is well known that a horizontal vector field on $M$ is projectable to $N$ if and only if it is invariant by the principal action, and that every vector field on $M$ can be locally written as $\sum_j f_j X_j$, where $f_j$ are smooth functions and $X_j$ are invariant vector fields.

One can define an almost complex structure $J$ on $M$ in the following way. If $X$ is an invariant vector field on $M$, set
\begin{equation}\label{E:complex}
	\theta(JX) = - \pi^{*}( \eta( \pi_{*}X ) ),
	\qquad
	\pi_{*}(JX) = \phi(\pi_{*}X) + \check{\theta}_X \xi,
\end{equation}
where $\check{\theta}_X$ is the only smooth local function on $N$ such that $\check{\theta}_X \circ \pi = \theta(X)$.
Then, the definition of $J$ can be extended to any vector field on $M$ by the above remark.

By~\cite{Ogawa}, the almost complex structure $J$ is integrable if and only if the almost contact structure on $N$ is normal and $\dd\theta(X,JY) + \dd\theta(JX,Y) = 0$ for all $X,Y\in TM$.
The latter condition reads as $J\dd\theta = \dd\theta$, or  as $\phi\cdot F = F$. In particular $\iota_{\xi} F = 0$.

If the manifold $N$ is equipped with a compatible metric $g$, there is a natural metric $h$ on $M$, given by
\[
	h := \pi^{*} g + \theta\otimes\theta,
\]
which is $J$-Hermitian.
By a simple computation, the fundamental form $\omega = h(J\cdot, \cdot)$ can be written as
\[
	\omega = - \pi^*\Phi - \pi^* \eta \wedge \theta.
\]

Now let $N$ be an almost 3-contact metric manifold of dimension $4m+3$.
Then we have an almost \hc{} structure $(J_1, J_2, J_3)$ on $M$, which is integrable if and only if the almost 3-contact structure is \hn{} and $\phi_i \cdot F = F$ for $i=1,2,3$ (see also \cite{BGM-hc}).
Note that the quaternionic dimension of $M$ is $n=m+1$. 

Considering a compatible metric $g$ on $N$, the natural metric $h$ is \hh.
The conditions for $h$ to be (strong) HKT have been determined in~\cite{FFUV1}.
\begin{thm}\cite[Theorem 7.1]{FFUV1}\label{T:HKT}
The \hh{} structure $(J_1, J_2, J_3, h)$ is HKT if and only if
$\phi_r\cdot\dd\Phi_r - \dd\eta_r\wedge\eta_r$ coincide for $r = 1,2,3$, and 
$\phi_r \cdot \iota_{\xi_r}\dd\Phi_r$ coincide for $r=1,2,3$.
Moreover, this HKT structure is strong if and only if, for every $r = 1,2,3$, we have
\[
\begin{aligned}
	\dd (\phi_r\cdot (\iota_{\xi_r}\dd\Phi_r) ) &= 0, \\
	\dd (\phi_r\cdot \dd\Phi_r - \dd\eta_r\wedge\eta_r)
	&=	\left( - \phi_r\cdot (\iota_{\xi_r}\dd\Phi_r) + F \right)\wedge F .
\end{aligned}
\]
\end{thm}

\section{HKT-Einstein metrics}\label{S:HKTE}
In this section we study the HKT-Einstein equation for the \hh{} structures on circle bundles described previously.
In the basic case, they can be HKT-Einstein only if the Chern--Ricci curvature vanishes.
However, we can consider suitable deformations of the \hh{} structure, that lead to examples with positive Einstein constant. 
Finally, we observe that these \hh{} manifolds always have constant Chern scalar curvature.

\subsection{The Lee form}\label{SS:Setup}
We begin by considering a normal, almost contact metric manifold $N$, of dimension $2m+1$, and a circle bundle $\pi \colon M \to N$ with the natural Hermitian structure $(J,h)$ as in Section~\ref{SS:cb}.
We compute the Lee form $\tau = - J\!\dd^{*}\omega$.
\begin{lem}\label{L:LeeForm}
\[
	\tau = (\pi^*b)\theta + \pi^*(\gamma - a \eta),
\]
where $a,b$ are smooth functions and $\gamma$ is a 1-form on $N$, given by
\[
	a := \frac{1}{2}\sum_{j=1}^{2m} F(\phi u_j , u_j),
	\quad
	b := \frac{1}{2}\sum_{j=1}^{2m} \dd\eta(\phi u_j , u_j),
	\quad
	\gamma(X) := \frac{1}{2}\sum_{j=1}^{2m} \dd\Phi(X, \phi u_j , u_j),
\]
for a local $g$-orthonormal frame $\{u_j\}_{j=1}^{2m}$ of the distribution $\mathcal{D}$. 
\end{lem}
\begin{proof}
It is well known that
\[
	2 (\nabla^{h}_X\omega)(Y,Z) =
	2h( (\nabla^{h}_X J)Y, Z ) = \dd\omega(X,Y,Z) - \dd\omega(X, JY, JZ),
\]
and it follows immediately that
\[
	\tau(X) = - \frac{1}{2} \sum_{j=1}^{2n} \dd\omega(X,Je_j, e_j),
\]
where $\{e_j\}_{j=1}^{2n}$ is a local $h$-orthonormal frame of $TM$.
Here $n = m + 1$ and we take $\{\hat{u}_1, \dots, \hat{u}_{2m}, v, \hat{\xi}\}$, where
$\{u_j\}$ is a local $g$-orthonormal frame of the distribution $\mathcal{D}$ and $v$ is a vertical vector field such that $\theta(v)=1$ (note $Jv=\hat{\xi})$.

We consider first $X$ to be a fundamental vector field (hence vertical).
Since $\dd\omega = \pi^*(\eta\wedge F - \dd\Phi) - \pi^*(\dd\eta)\wedge\theta$, we have
\[
	\dd\omega(X, Jv, v) = - (\pi^*(\dd\eta)\wedge\theta) (X, Jv, v) = 0,
\]
and for every vector field $u$ on $N$
\[
\begin{aligned}
	\dd\omega(X, J\hat{u}, \hat{u}) 
	&= - (\pi^*(\dd\eta)\wedge\theta) (X, J\hat{u}, \hat{u}) 
	= - \theta(X) \pi^*( \dd\eta (\pi_*J\hat{u}, \pi_*\hat{u}) )\\
	&= - \theta(X) \pi^*( \dd\eta (\phi u, u) ).
\end{aligned}
\]
Thus
\[
\begin{aligned}
	\tau(X)
	&= 	- \frac{1}{2} \sum_{j=1}^{2m} \dd\omega(X,J\hat{u}_j, \hat{u}_j)
			- \frac{1}{2} \dd\omega(X,J\hat{\xi}, \hat{\xi})
			- \frac{1}{2} \dd\omega(X,Jv, v)\\
	&=		\frac{1}{2} \theta(X) \sum_{j=1}^{2m} \pi^*( \dd\eta (\phi u_j, u_j) )
	= \theta(X) \pi^*b.
\end{aligned}
\]
Next, we assume $X$ to be a horizontal and invariant vector field.
As above, we compute
\[
\begin{aligned}
	\dd\omega(X, Jv, v) 
	&= - (\pi^*(\dd\eta)\wedge\theta) (X, Jv, v)
	= - \theta(v) \pi^*( \dd\eta(\pi_*X, \pi_*Jv) )\\
	&= - \theta(v) \pi^*( \dd\eta(\pi_*X, \xi) )	= 0,
\end{aligned}
\]
where we used that $\iota_{\xi}\dd\eta=0$.
On the other hand, for a vector $u \in \mathcal{D}$ we have
\[
\begin{aligned}
	\pi^*(\eta\wedge F)(X, J\hat{u}, \hat{u})
	&= \pi^*(\eta\wedge F (\pi_*X, \phi u, u)\\
	&= \pi^*(	\eta(\pi_*X)F(\phi u, u)
				- \eta(\phi u) F(\pi_*X, u)
				+\eta(u) F(\pi_* X, \phi u)		)\\
	&= \pi^*(	\eta(\pi_*X)F(\phi u, u)	)\\
	\pi^*\dd\Phi (X, J\hat{u}, \hat{u})
	&= \pi^*( \dd\Phi (\pi_*X, \phi u, u) )\\
	\pi^*(\dd\eta)\wedge\theta(X, J\hat{u}, \hat{u}) 
	&= 0
\end{aligned}
\]
and moreover, using that $\phi\xi = 0$ and $\iota_{\xi}\dd\eta=0$,
\[
\begin{aligned}
	\pi^*(\eta\wedge F)(X, J\hat{\xi}, \hat{\xi})
	&=	\pi^*(\eta\wedge F (\pi_*X, \phi\xi, \xi)) = 0\\
	\pi^*\dd\Phi (X, J\hat{\xi}, \hat{\xi})
	&=	\pi^*(\dd\Phi (\pi_*X, \phi\xi, \xi)) = 0\\
	\pi^*(\dd\eta)\wedge\theta(X, J\hat{\xi}, \hat{\xi}) 
	&= \pi^*( \eta(\xi) \dd\eta(\pi_*X, \xi) ) = 0.
\end{aligned}
\]
It follows that
\[
\begin{aligned}
	\tau(X)
	&=	-\frac{1}{2} \sum_{j=1}^{2m} 
		\dd\omega( X, J\hat{u_j}, \hat{u_j})\\
	&=	-\frac{1}{2} \sum_{j=1}^{2m} 
		-\pi^*( \dd\Phi(\pi_*X, \phi u_j, u_j)) + \theta(JX) \pi^*( F(\phi u_j, u_j) )\\
	&= \pi^*\gamma(X) - \pi^*(a\eta)(X).
\end{aligned}
\]
\end{proof}

We see that $(M,J,h)$ is balanced if and only if $b=0$ and $\gamma = a\eta$.
In particular
\begin{itemize}
\item
if $N$ is quasi-Sasakian, then $\gamma=0$, so $M$ is balanced if and only if $a=b=0$;
\item
if $N$ is $\alpha$-Sasakian, then $b=2m\alpha\neq0$, so $M$ is not balanced.
\end{itemize}

\subsection{The 3-quasi-Sasakian case}\label{SS:qS}
Let $N$ be a \hn{}, almost contact metric manifold, and let $\pi \colon M \to N$ be a circle bundle  with the natural \hh{} structure $(J_1, J_2, J_3, h)$ as above.
Let $4m+3$ be the dimension of $N$,  so that $n = m+1$ is the quaternionic dimension of $M$.
As specified in Section~\ref{SS:a3cm}, we assume $m\geq 1$, hence $n \geq 2$.

It is useful to note that $J_i\theta = - \pi^{*}\eta_i$ and $J_i(\pi^*\eta_j) = - \eps_{ijk}\pi^*\eta_k + \delta_{ij}\theta$ for every $i,j,k$.
We also have the following technical result.
\begin{lem}\label{L:J}
For every form $\psi$ on $N$, we have
\[
	J_r(\pi^{*}\psi) = \pi^{*}(\phi_r\cdot\psi) + \theta\wedge\pi^{*}(\phi_r \cdot \iota_{\xi_r}\psi).
\]
\end{lem}

\smallskip
Now, for every $X\in \Gamma(\calQ)$ we have $F(\phi_1 X, X) = F(\phi_2\phi_1 X, \phi_2X) = - F(\phi_1\phi_2X, \phi_2X)$.
Thus in Lemma~\ref{L:LeeForm} we obtain that $a=0$, and
the Lee form is  $\tau = (\pi^*b)\theta + \pi^*\gamma$,
where
\begin{equation}\label{E:b}
\begin{gathered}
	b 	:= \frac{1}{2}\sum_{j=1}^{4m} \dd\eta_1(\phi_1 u_j , u_j) 
		+ \frac{1}{2}\dd\eta_1(\phi_1 \xi_2 , \xi_2) 
		+ \frac{1}{2}\dd\eta_1(\phi_1 \xi_3 , \xi_3),	\\
	\gamma(X) 	:= \frac{1}{2}\sum_{j=1}^{4m} \dd\Phi_1(X, \phi_1 u_j , u_j)
					+ \frac{1}{2}\dd \Phi_1 (X, \phi_1 \xi_2 , \xi_2) 
					+ \frac{1}{2}\dd \Phi_1 (X, \phi_1 \xi_3 , \xi_3).
\end{gathered}
\end{equation}
Here $\{u_j\}$ is a local $g$-orthonormal frame of the distribution $\calQ$.\\

From now on, we shall assume that  $N$ is equipped with a 3-quasi-Sasakian structure.
For future use, we record the following:
\begin{lem}\label{L:rank}
Let $l$ be such that the rank of the 3-quasi-Sasakian structure on $N$ is $4l+1$ or $4l+3$.
Then $b = 2\delta(2l+1)$, and the Lee form of $h$ is $\tau=2\delta(2l+1)\,\theta$.
\end{lem}
\begin{proof}
As mentioned in Section~\ref{SS:a3cm}, if $\delta=0$ the structure on $N$ is 3-cosymplectic, hence $b=0$.
Therefore we consider $\delta\neq0$, and the rank to be $4l+3$.

Recall also that we have a splitting $\calQ = \calE^{4l}\oplus\calE^{4q}$, and that $\dd\eta_1 = 2\delta\Phi_1$ on $\calE^{4l}$.
The $g$-orthonormal local frame $\{u_j\}$ of $\calQ$ can be taken to have the first $4l$ elements in $\calE^{4l}$, and the remaining ones in $\calE^{4q}$. 
Thus
\[
	\sum_{j=1}^{4m} \dd\eta_1(\phi_1 u_j , u_j)  
	= \sum_{j=1}^{4l} \dd\eta_1(\phi_1 u_j , u_j) 
	= 2\delta\sum_{j=1}^{4l} g(\phi_1 u_j, \phi_1 u_j) 
	= 2\delta\sum_{j=1}^{4l} g(u_j, u_j) = 8l\delta.
\]
Since $\dd\eta_1(\phi_1\xi_2, \xi_2) + \dd\eta_1(\phi_1\xi_3, \xi_3) = 4\delta$, the claim is proved.
\end{proof}
Being $N$ 3-quasi-Sasakian, from Theorem~\ref{T:HKT} follows that $h$ is HKT if and only if 
\begin{equation}\label{E:eta}
\dd\eta_1 \wedge \eta_1 = \dd\eta_2 \wedge \eta_2 = \dd\eta_3 \wedge \eta_3.
\end{equation}
This condition is equivalent to $N$ being 3-$\delta$-cosymplectic.
Indeed, \eqref{E:eta} forces $\dd\eta_1$, $\dd\eta_2$, $\dd\eta_3$ to vanish when evaluated on sections of $\calQ$, giving
$\dd\eta_i = - 2\delta \eps_{ijk} \eta_j \wedge \eta_k$.
It is easy to check that the opposite implication holds as well.
See also~\cite[Corollary 4]{FFUV2}.

Therefore, assuming $h$ to be HKT, we have $b = 2\delta$.
The Lee form is then $\tau = 2\delta\theta$, and the first Chern-Ricci form of $(J_1, h)$ is
\[
	\RiC_{J_1} = \dd J_1 \tau = - 2\delta \pi^*(\dd\eta_1) = 4\delta^2 \pi^*(\eta_2\wedge\eta_3).
\]
Hence $J_2\RiC_{J_1} = - 4\delta^2 \theta \wedge \pi^* \eta_1$.
The HKT-Einstein equation takes the form
\[
	2\delta^2 \pi^*(\eta_2\wedge\eta_3) + 2 \delta^2 \theta\wedge\pi^*\eta_1
	= - \lambda \pi^*\Phi_1 + \lambda \theta\wedge\pi^*\eta_1.
\]
If this relation holds true, then evaluating it on a pair of vectors $X,Y\in\Gamma(\calQ)$, we find that necessarily $\lambda = 0$. This in turn implies $\delta = 0$, that is, $N$ is 3-cosymplectic. Note in this case $h$ is balanced and Chern--Ricci-flat.
Since the reversed implications hold as well, we get:
\begin{thm}\label{T:Einstein1}
Let $N$ be a 3-$\delta$-cosymplectic manifold, of dimension at least $7$, and let $M$ be a principal circle bundle over $N$ equipped with the HKT structure $(J_1, J_2, J_3, h)$ as in Section~\ref{SS:cb}.
The following facts are equivalent:
\begin{enumerate}
\item
$h$ is HKT-Einstein;
\item
$h$ is balanced, i.e., the Lee form vanishes;
\item
$h$ is Chern--Ricci-flat, i.e., $\RiC_{J_i}=0$ for all $i$;
\item
$N$ is 3-cosymplectic.
\end{enumerate}
\end{thm}
\begin{ex}
A well known example of a 3-cosymplectic manifold is the one given by Cabrera~\cite{Cabrera} on the torus $\mathrm{T}^{4m+3}$.
Considering a coframe $\{ e^j \}_{j=1}^{4m+3}$ whose elements are integral and closed, the structure equations are $\dd e^j = 0$ for every $j$.
As observed in \cite{Cabrera, FFUV2}, there exist non-trivial principal circle bundles over it, which correspond via Kobayashi's Theorem (\cite{Kob}) to $[F] \in \operatorname{H}^2(\mathrm{T}^{4m+3}, \Z)$, where
\[
	F = \sigma_{ijk} \sum_{l=1}^{m} 
		e^l \wedge e^{im+l} - e^{jm+l}\wedge e^{km+l}.
\] 
Here $\sigma_{ijk}$ denotes cyclic sum.
The total space $M$ can therefore be endowed with a balanced HKT structure.
\end{ex}

\subsection{Deformation of the \hh{} structure}\label{SS:mod}

Let us consider a manifold $N$ with a given 3-quasi-Sasakian structure $(\phi_i, \xi_i, \eta_i, g)$.
For a fixed $\mu\in\R^*$, we apply to it a \emph{$\calH$-homothetic deformation} in the sense of~\cite[Definition 2.3.1]{AgricolaDileo}, with parameters $a = 2\mu^2$, $b = - \mu^2$, $c = \mu$.
That is, we set
\[
	\tilde{g} = 2\mu^2 g  - \mu^2 \sum_{i=1}^{3}\eta_i\otimes\eta_i,
	\quad
	\tilde{\eta}_i = \mu \eta_i,
	\quad
	\tilde{\xi}_i = \frac{1}{\mu} \xi_i,
	\quad
	\tilde{\phi}_i = \phi_i.
\]
The metric is compatible with each almost contact structure.
Since the starting structure is hypernormal, from the definitions follows that also the deformed structure is hypernormal (of course, it will not be 3-quasi-Sasakian).

Consider a principal circle bundle $\pi \colon M \to N$ with a connection $\theta$, and as in the previous section let $(J_1, J_2, J_3, h)$ be the \hh{} structure on $M$ corresponding to $(\phi_i, \xi_i, \eta_i,g)$.
Now, if we perform the same construction, but starting from $(\tilde{\phi}_i, \tilde{\xi}_i, \tilde{\eta}_i, \tilde{g})$, we obtain a different \hh{} structure $(\tilde{J_1}, \tilde{J_2}, \tilde{J_3}, \tilde{h})$.
They are related as follows:
\[
\begin{gathered}
	\tilde{J}_i (\hat{u}) = J_i (\hat{u}) \text{ $\forall\,u\in\calQ$,}
	\qquad
	\tilde{J}_i (\hat{\xi}_i) = \mu J_i (\hat{\xi}_i),
	\qquad
	\tilde{J}_i (\hat{\xi}_j) = J_i (\hat{\xi_j}) \text{ if $j\neq i$,}\\
	\tilde{h} 
	= 2\mu^2 h - \mu^2\pi^{*}\!\left(\sum_{r=1}^3 \eta_r\otimes\eta_r \right) 
		+ (1-2\mu^2) \theta\otimes\theta.
\end{gathered}
\]
This specific deformation is chosen in order to obtain HKT structures.
\begin{lem}
The \hh{} structure $(\tilde{J_1}, \tilde{J_2}, \tilde{J_3}, \tilde{h})$ is HKT for every $\mu \in \R^*$.
\end{lem}
\begin{proof}
From~\cite{AgricolaDileo}, we know that the fundamental forms for $\tilde{g}$ are $\tilde{\Phi}_i = \mu^2 (2\Phi_i + \eta_j \wedge \eta_k)$ for every even permutation $(ijk)$ of $(123)$, thus
\begin{equation}\label{E:dPhi}
	\dd \tilde{\Phi}_i = \mu^2
		\left(\dd\eta_j \wedge \eta_k - \eta_j \wedge \dd\eta_k \right),
\end{equation}
Using~\eqref{E:formulas}, one obtains that 
\[
\begin{aligned}
	\tilde{\phi}_i \cdot \dd\tilde{\Phi}_i
	&=	\mu^2 (\phi_i \dd\eta_j) \wedge \eta_j +
		\mu^2 (\phi_i \dd\eta_k) \wedge \eta_k\\
	&= -\mu^2 \dd\eta_j \wedge \eta_j 
		- \mu^2 \dd\eta_k \wedge \eta_k 
		- 4\delta\mu^2 \,  \eta_i \wedge \eta_j \wedge \eta_k.
\end{aligned}
\]
Therefore we have the equality of $\tilde{\phi}_i \cdot \dd\tilde{\Phi}_i - \dd\tilde{\eta}_i \wedge \tilde{\eta}_i$ for $i = 1,2,3$.
It is  also easy to check that $\iota_{\tilde{\xi}_i} \dd\tilde{\Phi}_i = 0$ for every $i$, so we conclude applying Theorem~\ref{T:HKT}
\end{proof}
To determine the corresponding Lee form $\tilde{\tau}$, we use Lemma~\ref{L:LeeForm} and~\eqref{E:formulas}.
A straightforward computation yields
\[
	\tilde{\tau} = \frac{2\delta}{\mu}(l+1) \theta.
\]	
If $\delta = 0$,  the HKT structure $(\tilde{J_1}, \tilde{J_2}, \tilde{J_3}, \tilde{h})$ is balanced and hence Chern--Ricci-flat.
Therefore in the following we shall consider only the case $\delta \neq 0$.
\begin{thm}\label{T:3Sasakian}
Assume the 3-quasi-Sasakian structure on $N$ is not 3-cosymplectic (i.e., $\delta \neq 0$).
Then, the HKT structure $(\tilde{J_1}, \tilde{J_2}, \tilde{J_3}, \tilde{h})$ on $M$ is HKT-Einstein if and only if the 3-quasi-Sasakian structure on $N$ is 3-$\alpha$-Sasakian.
In this case, the Einstein constant is positive.
\end{thm}
\begin{proof}
We have
\[
\begin{aligned}
	\rho_{\tilde{J}_1}		
	&= \dd \tilde{J}_1 \tilde{\tau} 
		= - \frac{2\delta}{\mu}(l+1) \pi^*(\dd\tilde{\eta}_1)
		= - 2\delta(l+1) \pi^*(\dd\eta_1), \\
	\tilde{J}_{2}\rho_{\tilde{J}_1}	
	&= - 2\delta(l+1) \tilde{J}_2 \pi^*(\dd\eta_1)\\
	&=	- 2\delta(l+1)
		\left(
		\pi^*(\phi_2\cdot\dd\eta_1) 
		+ \theta\wedge \pi^*(\phi_2 \cdot \iota_{\tilde{\xi}_2}\dd\eta_1)
		\right)\\
	&=	- 2\delta(l+1)
		\left(
		- \pi^*\dd\eta_1 - 2\delta\pi^*(\eta_2\wedge\eta_3) +
		\frac{2\delta}{\mu}\theta\wedge(\pi^*\eta_1)
		\right),
\end{aligned}
\]
Here Lemma~\ref{L:J} was used, together with~\eqref{E:formulas}.
Taking into account that
\[
	\omega_{\tilde{J}_1}
	= 	- 2\mu^2 \pi^*\Phi_1 
		- \mu^2 \pi^*(\eta_2\wedge\eta_3)
		+ \mu \, \theta\wedge(\pi^*\eta_1),
\]
it is immediate to verify the HKT-Einstein equation, with a positive constant $\lambda > 0$, when $N$ is 3-$\alpha$-Sasakian ($\alpha\in\R^*$).

For the converse, let us assume $(\tilde{J_1}, \tilde{J_2}, \tilde{J_3}, \tilde{h})$ to be HKT-Einstein with Einstein factor $\lambda$.
Evaluating both sides of~\eqref{E:HKT-E} on the pair $(v, \hat{\xi}_1)$, where $v$ is a vertical vector field such that $\theta(v) = 1$, yields
\[
	\lambda = \frac{2\delta^2}{\mu^2} (l+1).
\]
Since we are assuming $\delta \neq 0$, it is $\lambda > 0$.
Inserting this value of $\lambda$, the equality reduces to $\dd\eta_1 = 2\delta\Phi_1$.
This in turn implies that the 3-quasi-Sasakian structure on $N$ is of maximal rank, hence it is 3-$\alpha$-Sasakian.
\end{proof}
This behavior reflects the fact that a 3-quasi-Sasakian manifold is locally the product of a 3-$\alpha$-Sasakian manifold (which is always Einstein) and a \hk{} manifold.

\begin{rmk}
There is an interesting relation with the class of \emph{3-$(\alpha,\delta)$-Sasakian manifolds} introduced in~\cite{AgricolaDileo}.
Indeed, applying the $\calH$-homothetic deformation of parameters $a=2\mu^2$, $b=-\mu^2$, $c=\mu$ to a 3-$\alpha$-Sasakian structure, 
the resulting structure is 3-$(\alpha', \delta')$-Sasakian, with $\alpha' = \alpha/2\mu$, $\delta' = \alpha/\mu$.
In other words, $\tilde{h}$ is the metric provided by the general construction in Section~\ref{SS:cb} when the base manifold $N$ is 3-$(\alpha, \delta)$-Sasakian, with $2\alpha = \delta$.

Conversely, starting from a 3-$(\alpha,\delta)$-Sasakian structure on $N$, it is easy to check using Theorem~\ref{T:HKT} that the \hh{} structure on $M$ is HKT only if  $2\alpha = \delta$.
(Note that, differently, the cone of a 3-$(\alpha, \delta)$-Sasakian manifold is HKT for every value of $\alpha$ and $\delta$ \cite[Corollary 4.3.1]{AgricolaDileo}).
\end{rmk}

\subsection{Examples of inhomogeneous HKT-Einstein manifolds}\label{SS:inhom}
We now use the last result to obtain examples of compact, non-homogeneous HKT-Einstein manifolds. 
We specialize to the cases in which $\mu = 1$ and the structure on $N$ is 3-Sasakian.
With this setup, we have
\[
	\tilde{J}_i = J_i, \quad i=1,2,3,	\quad\quad
	\tilde{h} 
	= 2 h - \pi^{*}\!\left(\sum_{r=1}^3 \eta_r\otimes\eta_r \right) 
		- \theta\otimes\theta.
\]
We remark that, in this particular case, the metric $\tilde{h}$ is the one appearing in~\cite[Theorem 5]{FFUV2}, up to homothety.

In~\cite{BGM-Main}, the authors develop a theory of reduction of 3-Sasakian manifolds, and apply it to build new simply connected examples, not necessarily homogeneous.
Let us briefly recall the construction appearing in \cite[Section 7]{BGM-Main}.
The sphere $\Sph^{4k-1}\subset \R^{4k} \simeq \HH^{k}$ with the canonical round metric $g_{\can}$ is a homogeneous 3-Sasakian manifold, with the group of 3-Sasakian automorphisms being $\Sp(k)$.
Let $T\subset\Sp(k)$ be the maximal torus consisting of diagonal matrices.
For every $p = (p_1, \dots, p_{k}) \in (\Z^*)^{k}$, the action of $T$ on $\Sph^{4k-1}$ restricts to the ``weighted'' circle action
\[
	\alpha_p \colon \Sph^{4k-1} \times \Sph^1 \to \Sph^{4k-1}
\]
given by $\alpha_p(u_1, \dots, u_k, t) = ( e^{\img 2\pi p_1 t}u_1, \dots,  e^{\img 2\pi p_n t}u_k) $, where $u_i$ are the quaternionic coordinates of $\HH^{k}$.

There exists a 3-Sasakian moment map $\mu_p$ for this action, and $M_p := \mu_p^{-1}(0)$ is a smooth compact submanifold of $\Sph^{4k-1}$ diffeomorphic to the complex Stiefel manifold $\SU(k)/\SU(k-2)$. 
Henceforth we assume $k\geq 3$, $0 < p_1 \leq\cdots\leq p_k$, and $p_1, \dots, p_k$ to be pairwise coprime.
In this way, the circle  action $\alpha_p$ restricts to a free action on $M_p$, and the quotient $N_p$ is a compact manifold.
Let $\iota_p \colon M_p \hookrightarrow \Sph^{4k-1}$ be the inclusion, and let $g_p$ be the unique metric on $N_p$ making the canonical projection $\pi\colon (M_p, \iota_p^*g_{\can}) \to (N_p, g_p)$ a Riemannian submersion.
The manifold $(N_p, g_p)$ turns out to be 3-Sasakian, and it is not a Riemannian homogeneous space, unless $p = (1,\dots, 1)$.

Consider the vector field $V_p$ generated by the free $\Sph^1$-action on $M_p$. 
Define a 1-form $\theta$  by decreeing that $\theta(V_p)=1$, and $\theta=0$ on the orthogonal complement of $V_p$ with respect to $\iota_p^*g_{\can}$.
This is a connection 1-form for the principal $\Sph^1$-bundle $\pi\colon M_p \to N_p$.
We have then the natural \hc{} structure on $M_p$, from $\theta$ and the 3-Sasakian structure on $N_p$.
In fact, $V_p$ is the restriction to $M_p$ of the vector field
\[
	\sum_{r=1}^k p_r	\left(
	u_r^{0}\frac{\del}{\del u_r^1} - u_r^{1}\frac{\del}{\del u_r^0} 
	+ u_r^{2}\frac{\del}{\del u_r^3} - u_r^{3}\frac{\del}{\del u_r^2}
	\right),
\]
where $u_r=u_r^0 + \img u_r^1 + \mathtt{j} u_r^2 + \mathtt{k} u_r^3$.
This makes possible to compute also $\theta$, and to verify the conditions on its curvature for the integrability of the \hc{} structure.

It is proved in~\cite{BGM-Stiefel} that this \hc{} structure is not homogeneous, for all $p\neq(1,\dots,1)$. 
Indeed, the fact that  $(N_p, g_p)$ is not homogeneous, implies that the metric $h$ on $M_p$ can not be invariant.
By \cite[Proposition 2.5]{BGM-Stiefel}, every hyper-complex automorphism of $M_p$ is an isometry for $h$,
 so the \hc{} structure on $M_p$ is not invariant as well.
Furthermore, for a fixed $k\geq3$, \hc{} structures corresponding to sequences $p$ and $q$ with the above constraints are equivalent if and only if $p=q$.

Therefore, as a consequence of Theorem~\ref{T:3Sasakian}, we have:
\begin{thm}\label{T:InHom}
For each $k\geq3$ and each $p$ with the conditions above, the non-homogeneous \hc{} manifolds $M_p$ admit a compatible  HKT-Einstein metric $\tilde{h}$ with positive Einstein factor.
\end{thm}
\begin{ex}\label{EX:SU(3)}
For $k=3$, the underlying manifold is the 8-dimensional Lie group $\SU(3)$.
The family of inhomogneous \hc{} structures is parametrized by a triple 
$(p_1, p_2, p_3)$ of pairwise coprime integers, with $0 < p_1 \leq p_2 \leq p_3$.
Each such structure admits a compatible HKT-Einstein metric with $\lambda=4$. 
\end{ex}
Recall that every \hc{} manifold has a unique torsion-free connection such that each complex structure is parallel, called \emph{Obata connection} (\cite{Obata}). 
By definition, its holonomy group  is contained in $\operatorname{GL}(n, \HH)$, where $n$ is the quaternionic dimension.
An important class of \hc{} manifolds is given by \emph{$\operatorname{SL}(n, \HH)$-manifolds}, those for which the Obata holonomy group lies in the commutator subgroup $\operatorname{SL}(n, \HH)$.
As observed in~\cite{FusiGentili}, HKT-Einstein metrics with positive factor can not be found on compact $\operatorname{SL}(n, \HH)$-manifolds.
Therefore:
\begin{cor}
The \hc{} manifolds $(M_p, J_1, J_2, J_3)$ are not $\operatorname{SL}(n, \HH)$.
\end{cor}
By contrast, when $p = (1,\dots, 1)$ we can recover the unique homogeneous HKT-Einstein metric on $\SU(k)/\SU(k-2)$ obtained in~\cite{BM}.

Indeed,  from~\cite[Section 6]{BGM-Main}, we know that  there exists a subgroup $G$ of the group of 3-Sasakian automorphisms of $\Sph^{4k-1}$, isomorphic to $\SU(k)$, acting transitively on $M_p$.
In this case, the action commutes with the circle action $\alpha_p$ and it descends to a transitive action on $N_p$, preserving $g_p$ and $\eta_1, \eta_2, \eta_3$.
Also the vector field $V_p$ on $M_p$ is invariant, and thus the connection 1-form $\theta$.
It follows that the metric $\tilde{h}$ is invariant.
The \hc{} structure is homogeneous as well, see~\cite[Theorem B]{BGM-Stiefel}, hence it is  a \emph{Joyce \hc{} structure}, as defined in~\cite{Joyce}.
Therefore we find a homogeneous HKT-Einstein metric $\tilde{h}$ on $\SU(k)/\SU(k-2)$ compatible with this \hc{} structure, which must necessarily coincide (up to homothety) with the unique one obtained in~\cite{BM}.
\begin{rmk}
If we consider a 3-Sasakian structure on $N$, and the trivial bundle $M = \Sph^1\times N$ with the canonical flat connection, the metric $h$ is known to be locally conformally \hk{} (LCHK) with Lee form parallel with respect to the Levi-Civita connection (\cite{BGM-Main, GOP}).
Up to scaling with a constant factor, $\tilde{h}$ is the HKT metric introduced in~\cite{OPS},  and~\cite[Proposition 2.1]{GentiliVezzoni} shows that $\tilde{h}$ is HKT-Einstein.
This fact can also be seen as a special case of Theorem~\ref{T:3Sasakian}.

We stress that, on the contrary, the bundles $M_p \to N_p$  considered in this section are not trivial; the manifolds $M_p$ are simply connected and do not admit compatible LCHK metrics (\cite{BGM-Stiefel}). 
Thus the HKT-Einstein metrics from Theorem~\ref{T:InHom} can not be recovered in that way.
\end{rmk}

\section{Quaternionic balanced and quaternionic Gauduchon metrics}\label{S:qBal}
We keep the same notations and conventions used so far, in particular $N$ is assumed to be 3-quasi-Sasakian, and we set $\dim N = 4m+3$, $m\geq1$, and $\dim M = 4m +4 = 4n$.
We saw in Section~\ref{SS:qS} that the \hh{} manifold $(M, J_1, J_2, J_3, h)$ is not always HKT. Thus, here we study the quaternionic balanced  and quaternionic Gauduchon conditions.\\

Let $\Omega$ be the associated $(2,0)$-form.
Our first goal is to compute the 1-form 
\[
	\beta_{\Omega} = \Lambda(\del\Omega)
\]
mentioned in  Section~\ref{SS:HHerm}, where $\Lambda$ is the adjoint of the Lefschetz operator $L= \Omega\wedge(-)$.
Consider a $h$-orthonormal local frame $\{e_1, \dots, e_{4n}\}$ adapted to $(J_1, J_2, J_3)$, and let $E_k= \frac{1}{\sqrt{2}}(e_{2k-1} - \img J_1 e_{2k-1}) = \frac{1}{\sqrt{2}}(e_{2k-1} - \img e_{2k})$ for $k=1,\dots,2n$, so that $J_2(E_{2j-1}) =\bar{E}_{2j}$.
Then $\Omega = \sum_{j=1}^{n} E^{2j-1}\wedge E^{2j}$. 
For a given $Z\in T_{J_1}^{1,0}M$ we have
\[
\begin{aligned}
	\beta_{\Omega}(Z)
	&= h( \iota_{Z}\!\del\Omega, \Omega)
	= \sum_{k=1}^{2n}\sum_{j<k} 
			\iota_{Z}\!\del\Omega(E_j, E_k) \Omega(E_j, E_k)\\
	&= \sum_{k=1}^{n}\dd\Omega(Z, E_{2k-1}, E_{2k}).
\end{aligned}
\]
(Recall that $\Omega$ is of type $(2,0)$ with respect to $J_1$.)
To compute the terms of the sum, we choose the basis $\{e_k\}_{k=1}^{4n}$ already defined in the proof of Proposition~\ref{P:ChernScal}.
Let $Z\in T_{J_1}^{1,0}M$ be such that its real and imaginary parts are projectable.
We have
\[
\begin{aligned}
	\Omega 
	&= \frac{1}{2}( \omega_{J_2} + \img\omega_{J_3} )
		= - \frac{1}{2}( 
		\pi^*\Phi_2 + \pi^*\eta_2\wedge\theta) 
		- \frac{\img}{2} (
		\pi^*\Phi_3 + \pi^*\eta_3\wedge\theta),\\
	\dd\Omega
	&=	\frac{1}{2}
		(\pi^*(\eta_2\wedge F) - \pi^*\dd\eta_2\wedge\theta ) 
		+ \frac{\img}{2} 
		(\pi^*(\eta_3\wedge F) - \pi^*\dd\eta_3\wedge\theta).
\end{aligned}	
\]
Consequently, the sum up to $m=n-1$ is
{\small
\[
\begin{aligned}
	&\sum_{k=1}^{m}\dd\Omega(Z, E_{2k-1}, E_{2k})\\
	&= \frac{1}{2}\sum_{k=1}^{m}
 			\dd\Omega(Z, e_{4k-3}, e_{4k-1}) 
 			- \dd\Omega(Z, e_{4k-2}, e_{4k})
 			- \img\! \dd\Omega(Z, e_{4k-2}, e_{4k-1}) 
 			- \img\! \dd\Omega(Z, e_{4k-3}, e_{4k})\\
	&= \frac{1}{4}  \left(\pi^*\eta_2(Z) + \img \pi^*\eta_3(Z) \right) \pi^*S_1
		-  \frac{1}{4}  \theta(Z) \pi^*S_2
		- \frac{\img}{4}  \theta(Z) \pi^*S_3,
\end{aligned}
\]
}
where
\[
\begin{aligned}
	S_1	
	&= \sum_{k=1}^m
			F(u_{4k-3}, u_{4k-1}) 
 			- F(u_{4k-2}, u_{4k})
 			- F(u_{4k-2}, u_{4k-1}) 
 			- F(u_{4k-3}, u_{4k})\\
	S_2
	&= \sum_{k=1}^m
			\dd\eta_2(u_{4k-3}, u_{4k-1}) 
 			- \dd\eta_2(u_{4k-2}, u_{4k})
 			+ \dd\eta_3(u_{4k-2}, u_{4k-1}) 
 			+ \dd\eta_3(u_{4k-3}, u_{4k})\\
	S_3
	&= \sum_{k=1}^m
			\dd\eta_3(u_{4k-3}, u_{4k-1}) 
 			- \dd\eta_3(u_{4k-2}, u_{4k})
 			- \dd\eta_2(u_{4k-2}, u_{4k-1}) 
 			- \dd\eta_2(u_{4k-3}, u_{4k})
\end{aligned}
\]
The fact that $F$ is $\phi_i$-invariant for every $i$ implies $S_1=0$.
Using~\eqref{E:formulas}, it is easy to see that
\[
	S_2 = - \sum_{k=1}^{m} \sum_{r=0}^{3}
			\dd\eta_1(\phi_1 (u_{4k-r}), u_{4k-r})
		= - 8\delta l,
	\qquad
	S_3 = 0,
\]
where $l$ is as in Lemma~\ref{L:rank}.
It follows that $\sum_{k=1}^{n-1}\dd\Omega(Z, E_{2k-1}, E_{2k}) = 2\delta l \, \theta(Z)$. 

Next, a straightforward calculation yields
\[
	\dd\Omega(X, E_{2n-1}, E_{2n})
	= 	 \delta\theta(X)
		+ \frac{\img}{4}\left(
		\pi^* (\dd\eta_3)(Z, \hat{\xi}_2) - \pi^*(\dd\eta_2)(Z, \hat{\xi}_3)
		\right).
\]
Recalling that $\iota_{\xi_i}\!\dd\eta_i = 0$ and that $\iota_{\xi_j}\dd\eta_i = -2\delta\eta_k$, we have
\[
\begin{aligned}
	\dd\eta_3(X, \xi_2) - \dd\eta_2(X, \xi_3)
	&= \eta_1(X)( \dd\eta_3(\xi_1, \xi_2) - \dd\eta_2(\xi_1, \xi_3) )\\
	&= \eta_1(X)( -\eta_3([\xi_1, \xi_2]) + \eta_2([\xi_1, \xi_3]) )\\
	&= - 4\delta\eta_1(X),
\end{aligned}	
\]
for every $X\in\Gamma(TN)$.
Thus
$\dd\Omega(Z, E_{2n-1}, E_{2n}) =  \delta\theta(Z)
- \img\delta \pi^*\eta_1(Z)$, and since $J_1 (\pi^*\eta_1) = \theta$, we obtain
\[
	\dd\Omega(Z, E_{2n-1}, E_{2n}) = 0.
\]
We conclude that 
\[
	\beta_{\Omega} = 2\delta l \,\theta^{1,0}.
\]
\begin{lem}\label{L:qbal}
The \hh{} structure $(J_1, J_2, J_3, h)$ is quaternionic balanced if and only if it is HKT.
\end{lem}
\begin{proof}
The quaternionic balanced condition is equivalent to $\beta_{\Omega}=0$, that is, $\delta=0$ or $l=0$. In both cases, the structure on $N$ is 3-$\delta$-cosymplectic, hence $h$ is HKT (see Section~\ref{SS:qS}).
\end{proof}
By the results in~\cite{FusiGentili}, the quaternionic Gauduchon condition is equivalent to
\begin{equation}\label{E:qbal}
	\scal^{\Bis} + 2\abs{\beta_{\Omega}}^2 = 0.
\end{equation}
Moreover, the Bismut-Ricci form associated to $(h, J_i)$ takes the form $\rho^{\Bis}_{J_i} = - \dd J_{i}(\beta_{\Omega} + \bar{\beta}_{\Omega})$, 
Here we have
\[
	\rho^{\Bis}_{J_i}
	= - 2\delta l \,\dd J_i \theta
	= 2\delta l \, \pi^*( \dd\eta_i ),
	\quad\quad
	\scal^{\Bis}= - 4\delta^2l(2l+1),
\]
where $\scal^{\Bis}$ can be determined as done for the Chern scalar curvature in Section~\ref{SS:ChernScal}.

We are ready to state the result of this section.
\begin{thm}
Let $N$ be a smooth manifold  with a 3-quasi-Sasakian structure, and let $M$ be a principal circle bundle over $N$ equipped with the \hh{} structure $(J_1, J_2, J_3, h)$ as in Section~\ref{SS:cb}.
The following three conditions are equivalent: $h$ is quaternionic Gauduchon, $h$ is quaternionic balanced, $h$ is HKT.
\end{thm}
\begin{proof}
We have to show that ~\eqref{E:qbal} holds if and only if $h$ is HKT.
Now $\abs{\theta^{1,0}}^2 = \sum_{k=1}^{2n} \abs{\theta^{1,0}(E_k)}^2$, where $E_k$ are defined as before.
The only non-zero term of the sum is $\theta^{1,0}(E_{2n-1})^2 = \frac{1}{2}$, so $\abs{\beta_{\Omega}}^2 = \frac{1}{2} \left( 2 \delta l \right)^2 = 2 \delta^2l^2$.
Using the above computation of $\scal^{\Bis}$ , one immediately notes that condition~\eqref{E:qbal} can be achieved only when $\delta=0$ or $l = 0$, so one can argue as in the proof of Lemma~\ref{L:qbal}.
\end{proof}

\subsection{Chern scalar curvature}\label{SS:ChernScal}
We are now able to explicitly determine the Chern scalar curvature of the \hh{} manifolds we are studying.
\begin{prop}\label{P:ChernScal}
The \hh{} manifold $(M, J_1, J_2, J_3, h)$ has constant non-negative Chern scalar curvature $\scal$, and $\scal=0$ if and only if $N$ is 3-cosymplectic.
The same holds for $(M, \tilde{J}_1, \tilde{J}_2, \tilde{J}_3, \tilde{h})$.
\end{prop}
\begin{proof}
We consider $(M, J_1, J_2, J_3, h)$. It holds that
\[
	\rho_{J_1}= \dd J_1( \tau - \beta_{\Omega} - \bar{\beta}_{\Omega} )
	= - 2\delta(l+1) \pi^*(\dd\eta_1).
\]
Choose a local unitary frame for $TM$ as follows.
Let $\{u_1, \dots, u_{4m}\}$ be a local $g$-orthonormal frame for $\calQ$ on $N$ such that
\[
\begin{aligned}
	\phi_1(u_{4k-3}) &= u_{4k-2}& 	\qquad 	\phi_1(u_{4k-1}) &= u_{4k} \\
	\phi_2(u_{4k-3}) &= u_{4k-1}& 	\qquad 	\phi_2(u_{4k-2}) &= - u_{4k} \\
	\phi_3(u_{4k-3}) &= u_{4k}&		\qquad 	\phi_3(u_{4k-2}) &= u_{4k-1}
\end{aligned}
\]
and let $e_j = \hat{u}_j$ for $j=1,\dots, 4m$ (recall $m = n-1$).
Set also $e_{4n-3}=v$, $e_{4n-2}=\hat{\xi}_1$, $e_{4n-1} = \hat{\xi}_2$, $e_{4n}=\hat{\xi}_3$, where again $v$ is the fundamental (vertical) vector field such that $\theta(v) = 1$.
As in the proof of Lemma~\ref{L:rank}, we find
\[
	\scal 
	= \frac{1}{2}\sum_{k = 1}^{4n} \rho_{J_1}(e_k, J_1 e_k)
	= 4\delta^2(l+1)(2l+1).
\]
The statement for $(M, \tilde{J}_1, \tilde{J}_2, \tilde{J}_3, \tilde{h})$ can be proved similarly.
\end{proof}
We saw in Section~\ref{S:HKTE} that there are cases in which the structures $(J_1, J_2, J_3, h)$, as well as their deformations $(\tilde{J}_1, \tilde{J}_2, \tilde{J}_3, \tilde{h})$, are HKT but not HKT-Einstein. 
Nevertheless, this proposition ensures they lie in the strictly larger family of HKT structures with constant Chern scalar curvature.

\section{HKT-Einstein metrics and the twist construction}\label{S:Twist}
Here we study the HKT-Einstein equation on manifolds obtained via \emph{twist construction}.
The latter was systematically investigated by A. Swann in~\cite{Swann}, but see also~\cite{GrantcharovPoon, Joyce}.

The construction can be summarized as follows\footnote{The notation here differs from that in~\cite{Swann} to avoid ambiguity with the previous sections.}.
Given a smooth manifold $M$ with a $\Sph^1$-action, consider a principal $\Sph^1$-bundle $\pi \colon P \to M$ with a connection.
Let $\pi^*\Psi$ be its curvature, where $\Psi$ is a closed integral $2$-form on $M$.
Denoting by $X$ a fundamental vector field generating the action on $M$,
suppose there exists $f\in\operatorname{C}^{\infty}(M;\R)$ such that $\dd f = - \iota_{X}\Psi$.
In particular $\mathcal{L}_X\Psi = 0$. 
We assume $f$ to be nowhere-vanishing; in the compact case this can be done  with no loss of generality.
As proved in~\cite[Proposition 2.1]{Swann}, the $\Sph^1$-action on $M$ can be lifted to one on $P$, commuting with the principal action and preserving the connection.
We assume this action to be free, so that the quotient $M':=P/\Sph^1$ is a smooth manifold.
There is another principal circle bundle $\pi'\colon P \to M'$. 

Let $\calH$ be the horizontal distribution on $P$ given by the connection, and let the lifted action be transverse to $\calH$.
An invariant $(p,0)$-tensor $\alpha$ on $M$ determines uniquely a tensor $\alpha'$ on $M'$ such that $\pi^*\alpha = {\pi'}^* \alpha'$ on $\calH$. 
One says that $\alpha$ and $\alpha'$ are \emph{$\calH$-related}, in symbols $\alpha\sim_{\calH}\alpha'$.
In case $\alpha$ is a differential form, by~\cite[Corollary 3.6]{Swann} the exterior derivatives are related as
\begin{equation}\label{E:drelation}
	\dd\alpha' \sim_{\calH} 
	\dd\alpha - f^{-1} \Psi\wedge \iota_X\alpha.
\end{equation}
If $M$ is equipped with an invariant complex structure $J$, there is a corresponding almost complex structure $J'$ on $M'$, see~\cite[Section 3]{Swann} for the definition. In order to have $J'$  integrable, we consider $\Psi$ to be of type $(1,1)$.
Also, if $J$ is compatible with a metric $h$ on $M$, then $J'$ is compatible with $h'$  and the respective fundamental forms $\omega$ and $\omega'$ are $\calH$-related.

We can find the relation between the respective Lee forms $\tau$ and $\tau'$.
As in the proof of Lemma~\ref{L:LeeForm}, we have
\[
	\tau'(V) = - \frac{1}{2} \sum_{j=1}^{2n} \dd\omega'(V, J'e'_j, e'_j),
	\quad V\in\Gamma(TM'),
\]
where $\{e'_j\}_{j=1}^{2n}$ is a local $h'$-orthonormal frame of $TM'$ and $n = \dim_{\C}M = \dim_{\C}M'$.
Let $U$ be the vector field on $M$ determined by the property that its horizontal lift to $\calH$ coincides with the horizontal lift of $V$.
Let $\{e_j\}$ be the $h$-orthonormal local frame of $TM$ obtained from $\{e'_j\}$ in the same way.
Note that by definition of $J'$ the horizontal lift of $J'e'_j$ coincides with the horizontal lift of $Je_j$.
Therefore we have
\[
\begin{aligned}
	{\pi'}^* (\tau'(V))
	&=  - \frac{1}{2} \sum_{j=1}^{2n} 
			{\pi'}^*(\dd\omega'(V, J'e'_j, e'_j) )\\
	&= - \frac{1}{2} \sum_{j=1}^{2n} 
		\pi^*
		\left(
		(\dd\omega - f^{-1}\Psi\wedge\iota_X\omega)
		(U, Je_j, e_j)
		\right)\\
	&= \pi^*( \tau(U) ) 
		+ \frac{1}{2f} \sum_{j=1}^{2n} 
		\pi^*
		\left(
		\Psi\wedge\iota_X\omega (U, Je_j, e_j)
		\right).
\end{aligned}
\]
We compute
\[
\begin{aligned}
	(\Psi\wedge\iota_X\omega)(U, Je_j, e_j)
	&= \Psi(U, Je_j)\omega(X, e_j)
		+ \Psi(Je_j, e_j)\omega(X, U)
		+ \Psi(e_j, U)\omega(X, Je_j)\\
	&= \Psi(Je_j, U) h(Je_j, X)
		+ \Psi(e_j, U) h(X, e_j)
		+ \Psi(e_j, Je_j) h(JU, X),
\end{aligned}
\]
and it follows that
\[
	{\pi'}^* \tau'
	= \pi^* 
		\left( 	\tau + \frac{1}{f} \iota_{X}\Psi +
		\frac{1}{2f} ( \operatorname{tr}_{\omega}\Psi ) J(X^{\flat})
		\right),
\]
where $\operatorname{tr}_{\omega}\Psi = \sum_{j}\Psi(e_j, Je_j)$.

Now let $M$ be equipped with an invariant HKT structure $(J_1, J_2, J_3, h)$.
We assume $\Psi$ to be of type $(1,1)$ with respect to each complex structure $J_1, J_2, J_3$.
The induced \hh{} structure on the twist $M'$ is HKT, by~\cite[Proposition 5.2]{Swann}.
Moreover, it is immediate to see that $\tr_{\omega_{J_1}}\Psi = 0$.

In order to study the HKT-Einstein equation, we compute the first Chern-Ricci form of the Hermitian structure $(J'_1, h')$, that we denote simply by $\RiC_{J'_1}$.
We have 
\[
	J'_1\tau' \sim_{\calH} J_1\tau + \frac{1}{f} J_1 (\iota_X\Psi) 
	= J_1\tau - \frac{1}{f} \iota_{J_1X}\Psi,
\] 
hence by \eqref{E:drelation}
\[
\begin{aligned}
	\RiC_{J'_1} 
	&= \dd J'_1\tau'\\
	&\sim_{\calH} \dd\left( J_1\tau - \frac{1}{f}\iota_{J_1X}\Psi \right)
		- \frac{1}{f} \Psi \wedge \iota_{X}
			\left( J_1\tau - \frac{1}{f}\iota_{J_1X}\Psi \right)\\
	&=	\RiC_{J_1} 
		- \frac{1}{f^2} (\iota_X\Psi)\wedge(\iota_{J_1X}\Psi)
		- \frac{1}{f} \dd\iota_{J_1X}\Psi
		- \frac{1}{f} \tau(J_1X) \Psi
		+ \frac{1}{f^2} (\iota_{X} \iota_{J_1X}\Psi) \Psi.
\end{aligned}
\]
Taking into account that 
\[
	J_2(\iota_{X}\Psi) = - \iota_{J_2X}\Psi,
	\quad
	J_2(\iota_{J_1X}\Psi) = \iota_{J_3X}\Psi, 
	\quad\text{and}\quad
	\dd\iota_{J_1X}\Psi = \mathcal{L}_{J_1X}\Psi,
\]
we obtain
\[
\begin{aligned}
	\RiC_{J'_1} - J'_2 \RiC_{J'_1}
	&= \left( \RiC_{J_1} - J_2\RiC_{J_1} \right)
	- \frac{1}{f^2}
	(\iota_X\Psi)\wedge(\iota_{J_1X}\Psi)
	- \frac{1}{f^2} (\iota_{J_2X}\Psi)\wedge(\iota_{J_3X}\Psi)\\
	&\quad - \frac{1}{f} 
	( \mathcal{L}_{J_1X}\Psi - J_2 \mathcal{L}_{J_1X}\Psi ).
\end{aligned}
\]
The following result is now an immediate consequence of this computation.
\begin{thm}\label{T:EinsteinTwist}
Let $M$ be a HKT manifold, and let $M'$ be the twist associated to the twist data $(X, \Psi, f)$, where the $\Sph^1$-action generated by $X$ preserves the HKT structure, and $\Psi$ is of type $(1,1)$ with respect to each complex structure.
Suppose in addition that
\[
	(\iota_X\Psi)\wedge(\iota_{J_1X}\Psi)
	+ (\iota_{J_2X}\Psi)\wedge(\iota_{J_3X}\Psi)
	= - f (\mathcal{L}_{J_1X}\Psi - J_2 \mathcal{L}_{J_1X}\Psi).
\]
Then the HKT structure on $M'$ is HKT-Einstein if and only if the corresponding structure on $M$ is so.
\end{thm}
\begin{cor}\label{C:EinsteinTwist}
Let $M$ be a HKT-Einstein manifold.
Consider the twist data $(X, \Psi, f)$, where the $\Sph^1$-action generated by $X$ preserves the HKT structure, $\Psi$ is of type $(1,1)$ with respect to each complex structure, and $\iota_{X}\Psi = 0$, that is, $f$ is constant.
Then the twist $M'$ is HKT-Einstein (with same Einstein factor).
\end{cor}
\begin{proof}
We just  note that $\iota_X\Psi=0$ implies the assumption in Theorem~\ref{T:EinsteinTwist}.
In fact, it holds that $\tau'\sim_{\calH}\tau$ and $(\RiC_{J'_1} - J'_2\RiC_{J'_1}) \sim_{\calH}(\RiC_{J_1} - J_2\RiC_{J_1})$.
\end{proof}
\medskip
Let us compare the twist construction with the class of manifolds examined in the previous sections.

We consider a \hh{} manifold $(M, J_1, J_2, J_3, h)$ obtained from a principal $\Sph^1$-bundle $\pi_M \colon M \to N$  with connection $\theta_M$, and from a 3-quasi-Sasakian structure $(\phi_i, \xi_i, \eta_i, g)$ on $N$.
Let $\pi_{P_\circ} \colon P_{\circ} \to N$ be another principal circle bundle, with a connection $\theta_{P_{\circ}}$.
We call $F$ and  $\Psi_{\circ}$ the curvature forms respectively of $\theta_M$ and $\theta_{P_{\circ}}$, and we assume them to be invariant under the endomorphisms $\phi_i$, $i=1,2,3$.

Let $P$ be the pull-back bundle of $P_{\circ}$ under the map $\pi_{M}$, namely
\[
	P := \{ (x,p) \in M\times P_{\circ} 
			\colon \pi_M(x) = \pi_{P_{\circ}}(p) \},
\]
with the map $\pi_P \colon P \to M$ being the canonical projection on the first factor.
This bundle can be equipped with the connection $\theta_P :=  \psi^*(\theta_{P_{\circ}})$, where $\psi$ is the projection on the second factor.
Its curvature is $\Psi := \pi_M^*\Psi_{\circ}$, which by construction is of type $(1,1)$ for each $J_i$.
Since $M$ is the total space of a principal $\Sph^1$-bundle, it admits a free $\Sph^1$-action, preserving the \hh{} structure.
Hence on $P$ we have the lifted free circle action 
\[
	z\cdot(x,p) = (z\cdot x, z^{-1}\cdot p),
\]
which covers the action on $M$, preserves $\theta_P$, and commutes with the principal action on $P$.

Let $(M', J'_1, J'_2, J'_3, h')$ be the twist of $M$ via the above data, and let 
$\pi'_{P} \colon P \to M'$ be the projection.
Like $M$, also $M'$ has the structure of a principal $\Sph^1$-bundle over $N$. To see this explicitly, write points of $M'$ as $[x, p]$ (orbits of the lifted $\Sph^1$-action on $P$), and put $\pi'_{M} \colon M' \to N$, $\pi'_{M}([x,p]) = \pi_M(x)$.
Now, the unique 1-form on $M'$ $\calH$-related to $\theta_M$ defines a connection for this bundle; we call it $\theta_{M'}$.
As $\pi_{M}\circ\pi_{P} = \pi'_{M}\circ\pi'_{P}$, one immediately proves that $h' = (\pi'_{M})^* g + \theta_{M'}\otimes\theta_{M'}$.
Similarly, the complex structures $J'_i$ obey the relations in~\eqref{E:complex} with projection $\pi_{M'}$  and  connection $\theta_{M'}$.
Altogether, this shows the following fact.
\begin{prop}
The twisted manifold $(M', J'_1, J'_2, J'_3, h')$ can be equivalently obtained via the method in Section~\ref{SS:cb}, considering the bundle $\pi'_{M}\colon M' \to N$ and the connection $\theta_{M'}$.
\end{prop}
We can also view \hh{} structures on non-trivial circle bundles as twists of trivial circle bundles.
\begin{prop}\label{P:trivial}
The \hh{} manifold $(M, J_1, J_2, J_3, h)$ is the twist of the manifold $\Sph^1\times N$ with the \hh{} structure associated to the canonical flat connection.
\end{prop}
\begin{proof}
For $M$ fixed, consider the bundle $P_{\circ} = M \to N$. Its pullback to $\Sph^1\times N$ is $P = \Sph^1\times M$. 
The above argument proves that twisting $\Sph^1\times N$ with the \hh{} structure associated to the canonical connection, via the bundle $P$, gives back the initial $M$.
\end{proof}
\begin{rmk}
If we start with a HKT-Einstein manifold $(M, \tilde{J}_1, \tilde{J}_2, \tilde{J}_3, \tilde{h})$ coming from a 3-$\alpha$-Sasakian manifold $N$, as in Theorem~\ref{T:3Sasakian}, the twist will be again a HKT-Einstein manifold from the same family.
Indeed, $\iota_X\Psi=0$ by construction , and we can apply Corollary~\ref{C:EinsteinTwist}. Note that in this particular case we also have $\rho_{J'_1}\sim_{\calH}\rho_{J_1}$, since $\tau(J_1X)=0$.

On the other hand, combining Corollary~\ref{C:EinsteinTwist} with Proposition~\ref{P:trivial}, we see that the HKT-Einstein manifold $(M, \tilde{J}_1, \tilde{J}_2, \tilde{J}_3, \tilde{h})$ can be recovered from the canonical HKT-Einstein structure on the trivial circle bundle $\Sph^1\times N$ with a suitable twist. This provides an alternative proof of part of Theorem~\ref{T:3Sasakian}.
\end{rmk}
More examples of positive HKT-Einstein manifolds can be obtained twisting manifolds of the form $M = M_1 \times M_2$, where $M_1$ is a HKT-Einstein manifold of the type of Theorem~\ref{T:3Sasakian}, while $M_2$ is a Joyce hypercomplex manifold with its unique invariant HKT-Einstein metric (\cite{BM}).
The metrics on the two factors can be scaled to have equal Einstein factors, so that the product metric on $M$ is HKT-Einstein.
In order to apply Corollary~\ref{C:EinsteinTwist}, the twist data can be chosen as follows.

$\bullet$
For the 2-form $\Psi$, we may consider the pull-back to $M$ of a closed 2-form $\Psi_1$ on $M_1$ invariant by each complex structure and representing an integral cohomology class.
As noted by Swann~\cite{Swann}, these can be found when $M_1$ is a circle bundle over a 3-Sasakian manifold $N$, since harmonic 2-forms on $N$ are invariant under $\phi_i$, for every $i$ (\cite{GalickiSalamon}).
The twisting bundle is then $P=P_1\times M_2 \to M = M_1\times M_2$, where
$P_1 \to M_1$ is a principal circle bundle corresponding to $[\Psi_1]$.

$\bullet$
Furthermore, $M_2$ admits invariant vector fields that are hyper-holomorphic and Killing for the invariant HKT-Einstein metric.
For simplicity, let us consider the case of a Lie group of the form $T^{\ell}\times G$, where $G$ is a semisimple compact Lie group.
Then, the existence of hyper-holomorphic vector fields is due to~\cite[Lemma 3.2]{BrienzaGentiliEtAl}. 
It is possible to show that they are also Killing for the invariant HKT-Einstein metric thanks to the explicit structure of this metric, see~\cite{BM}.
This fact can be generalized to some coset spaces.

In any case, let $X_2$ be one of such vector fields, and let $X$ be the field induced on $M$.
The $\Sph^1$-action on $M$ generated by $X$ is given by the action on $M_2$ generated by $X_2$, extended to act trivially on the factor $M_1$.

\begin{ex}
Here is a low-dimensional example.
Let $M_1$ be as in Example~\ref{EX:SU(3)}, for a  triple $p=(p_1,p_2,p_3)\neq(1,1,1)$. 
Then $M_1$ is diffeomorphic to $\SU(3)$ and has a inhomogeneous \hc{} structure; here $P_1 = \Sph^1\times\SU(3)$.
For $M_2$, we may consider again $\SU(3)$, but endow it with a left-invariant Joyce \hc{} structure. 
The latter is associated to a decomposition $\fraksu(3)=\R X\oplus \frakd_1\oplus\frakf_1$, with notations as in~\cite{BM}. Then $X$ gives rise to a left-invariant vector field, that is hyper-holomorphic and Killing for the compatible left-invariant HKT-Einstein metric. It generates a circle subgroup, denoted by $\U(1)$.
Therefore, we have inhomogeneous HKT-Einstein structures on
\[
	\Sph^1 \times \SU(3) \times \frac{\SU(3)}{\U(1)}.
\]
\end{ex}

\end{document}